\documentclass{gtmon_a}
\pdfoutput=1
\usepackage[all]{xy}

\proceedingstitle{Groups, homotopy and configuration spaces (Tokyo
  2005)}
\conferencestart{5 July 2005}
\conferenceend{11 July 2005}
\conferencename{Groups, homotopy and configuration spaces,
                in honour of Fred Cohen's 60th birthday}
\conferencelocation{University of Tokyo, Japan}

\editor{Norio Iwase}
\givenname{Norio}
\surname{Iwase}

\editor{Toshitake Kohno}
\givenname{Toshitake}
\surname{Kohno}

\editor{Ran Levi}
\givenname{Ran}
\surname{Levi}

\editor{Dai Tamaki}
\givenname{Dai}
\surname{Tamaki}

\editor{Jie Wu}
\givenname{Jie}
\surname{Wu}

\title[Euclidean configuration spaces and associated fibrations]{On 
       the category of Euclidean configuration spaces\\and associated 
       fibrations}
\author{Fridolin Roth}
\givenname{Fridolin}
\surname{Roth}
\address{Fachbereich Mathematik der Universit\"at Hamburg\\
SPAZ\\Bundesstrasse 55\\\newline
20146 Hamburg\\Germany}
\email{roth@math.uni-hamburg.de}
\urladdr{http://www.math.uni-hamburg.de/home/roth/}

\volumenumber{13}
\issuenumber{}
\publicationyear{2008}
\papernumber{20}
\startpage{447}
\endpage{461}

\doi{}
\MR{}
\Zbl{}

\arxivreference{}

\keyword{Schwarz genus}
\keyword{sectional category}
\keyword{Lusternik--Schnirelmann category}
\keyword{category weight}
\keyword{configuration spaces}

\subject{primary}{msc2000}{55M30} 
\subject{secondary}{msc2000}{55R80}
\subject{secondary}{msc2000}{55S40}

\received{31 May 2006}
\revised{22 March 2007}
\accepted{29 March 2007}
\published{19 March 2008}
\publishedonline{19 March 2008}
\proposed{}
\seconded{}
\corresponding{}
\version{}

\makeatletter
\def\cnewtheorem#1[#2]#3{\newtheorem{#1}{#3}[section]
\expandafter\let\csname c@#1\endcsname\c@thm}

\let\xysavmatrix\xymatrix
\def\xymatrix{\disablesubscriptcorrection\xysavmatrix}
\AtBeginDocument{\let\bar\wbar\let\tilde\wtilde\let\hat\what}

\newtheorem{thm}{Theorem}[section]
\cnewtheorem{conj}[thm]{Conjecture}
\cnewtheorem{prop}[thm]{Proposition}
\cnewtheorem{lem}[thm]{Lemma}
\cnewtheorem{cor}[thm]{Corollary}

\theoremstyle{definition}
\cnewtheorem{Def}[thm]{Definition}

\makeatother  

\newcommand{\sym}{\Sigma}

\renewcommand{\R}{\mathbb R}
\newcommand{\N}{\mathbb N}
\renewcommand{\Z}{\mathbb Z}

\newcommand{\col}{\colon\thinspace}

\DeclareMathOperator{\cat}{cat}
\DeclareMathOperator{\catcl}{cat^{\mathit{cl}}}

\DeclareMathOperator{\secat}{secat}

\DeclareMathOperator{\wgt}{wgt}
\DeclareMathOperator{\cuplength}{cup}

\newcommand{\zz}{\mathbb Z / 2 \mathbb Z}
\newcommand{\zp}{\mathbb Z / p \mathbb Z}

\newenvironment{fat}{\begin{bfseries}}{\end{bfseries}}

\DeclareMathOperator{\cohdim}{cohdim}

\DeclareMathOperator{\tth}{th}
\DeclareMathOperator{\tst}{st}

\begin{document}

\begin{htmlabstract}
We calculate the Lusternik&ndash;Schnirelmann category of the
k<sup>th</sup> ordered configuration spaces F(<b>R</b><sup>n</sup>,k)
of <b>R</b><sup>n</sup> and give bounds for the category of the
corresponding unordered configuration spaces B(<b>R</b><sup>n</sup>,k)
and the sectional category of the fibrations
&pi;<sup>n</sup><sub>k</sub>:F(<b>R</b><sup>n</sup>,k)&rarr;
B(<b>R</b><sup>n</sup>,k). We show that
secat(&pi;<sup>n</sup><sub>k</sub>) can be expressed in terms of
subspace category. In many cases, eg, if n is a power of 2, we
determine cat(B(<b>R</b><sup>n</sup>,k)) and
secat(&pi;<sup>n</sup><sub>k</sub>) precisely.
\end{htmlabstract}

\begin{abstract}
We calculate the Lusternik--Schnirelmann category of the $k^{\rm th}$
ordered configuration spaces $F(\mathbb{R}^n,k)$ of $\mathbb{R}^n$ and
give bounds for the category of the corresponding unordered
configuration spaces $B(\mathbb{R}^n,k)$ and the sectional category of
the fibrations $\pi^n_k\colon\thinspace F(\mathbb{R}^n,k)\rightarrow
B(\mathbb{R}^n,k)$. We show that ${\rm secat}(\pi^n_k)$ can be
expressed in terms of subspace category. In many cases, eg, if $n$ is
a power of $2$, we determine ${\rm cat}(B(\mathbb{R}^n,k))$ and ${\rm
secat}(\pi^n_k)$ precisely.
\end{abstract}

\begin{asciiabstract}
We calculate the Lusternik-Schnirelmann category of the k-th ordered
configuration spaces F(R^n,k) of R^n and give bounds for the category
of the corresponding unordered configuration spaces B(R^n,k) and the
sectional category of the fibrations pi^n_k: F(R^n,k) --> B(R^n,k). We
show that secat(pi^n_k) can be expressed in terms of subspace
category. In many cases, eg, if n is a power of 2, we determine
cat(B(R^n,k)) and secat(pi^n_k) precisely.
\end{asciiabstract}

\maketitle

\section{Motivation and results}

The sectional category $\secat(p)$ of a fibration $p \col E\rightarrow
B$ is defined to be the least integer $n$ such that the base $B$ can
be covered with $n+1$ open sets admitting local sections (Cornea,
Lupton, Oprea, and Tanr{\'e} \cite[9.13]{CLOT03}). This notion, as in
{\v{S}}varc \cite[page 70]{Sva} sometimes defined without the
$+1$--shift and referred to as the Schwarz genus of $p$, has proved
very useful. For instance, consider the fibration $\pi^2_k\col
F(\R^2,k)\rightarrow B(\R^2,k)$ from the ordered configuration space
$F(\R^2,k)=\{(x_1,x_2,\dots,x_k)\in (\R^2)^k| x_i\neq x_j \mbox{ for }
i\neq j\}$ to its unordered quotient $B(\R^2,k):=F(\R^2,k)/ \sym_k$
given by the obvious action of the symmetric group $\sym_k$ on $k$
letters.  The sectional category of the fibration $\pi^2_k\col
F(\R^2,k)\rightarrow B(\R^2,k)$ has attracted some attention since
Smale showed that it provides a lower bound for the complexity of
algorithms computing the (pairwise disjoint) roots of a complex
polynomial of degree $k$ (Smale \cite{Sma87}, Vasiliev
\cite{Vas88,Vas92}, de Concini, Procesi and Salvetti \cite{DeCPS04}
and Arone \cite{Aro05}). Nevertheless $\secat(\pi^2_k)$ has not yet
been determined for all $k$.

We now recollect what is known on $\secat(\pi^2_k)$ so far: Vassiliev
showed the inequality $k-D_p(k)\leq \secat(\pi^2_k)\leq k-1$ where
$D_p(k)$ is the sum of the coefficients in the $p$--adic extension of
$k$ \cite{Vas88,Vas92}. In particular, this gives $\secat(\pi^2_k)=
k-1$ if $k$ is a power of a prime. If $k$ is not a power of a prime,
it is very difficult to determine $\secat(\pi^2_k)$ precisely, and it
is only recently that some progress has been made. In \cite{DeCPS04}
de Concini, Procesi and Salvetti developed an obstruction theory to
decide whether $\secat(\pi^2_k)$ equals its known upper bound $k-1$
and showed that for $k=6$ --- the least $k$ for which the question was
open --- in fact it does not. Based on that theory, Gregory Arone did
some more calculations and showed that $\secat(\pi^2_k)<k-1$ holds for
all $k$ which are neither the power of a prime nor twice the power of
a prime \cite{Aro05}. If $k=2p^{\ell}$, the question whether
$\secat(\pi^2_k)=k-1$ is still open for some $\ell$ and odd $p$, as
well as the precise determination of $\secat(\pi^2_k)$ in many other
cases.

In this paper we begin to study the sectional category
$\secat(\pi^n_k)$ of the fibrations $\pi^n_k\col F(\R^n,k)\rightarrow
B(\R^n,k)$ for varying $n\in \N$. This is closely related to the
Lusternik--Schnirelman category $\cat(B(\R^n,k))$ of the unordered
configuration spaces $B(\R^n,k)$. Here we follow \cite{CLOT03} and say
that the Lusternik--Schnirelman category $\cat(X)$ of a topological
space $X$ is the least integer $m$ such that $X$ can be covered with
$m+1$ open sets, which are all contractible within $X$.  One
elementary relation between $\secat(\pi^n_k)$ and $\cat(B(\R^n,k))$ is
the general fact that the sectional category of a fibration is bounded
above by the category of its base. Together with
$\cat(B(\R^n,k))\leq(k-1)\cdot(n-1)$ (\fullref{upperbound}) this
gives an upper bound for $\secat(\pi^n_k)$. Moreover, in our cases, we
get descriptions of $\secat(\pi^n_k)$ in terms of the category of a
map and subspace category, definitions whereof are given in \fullref{LSC}:

\begin{thm} \label{sec+sub}
Let $n,k,r \in \{1,2,3,\dots\}$ 
and let $\pi^n_k\col  F(\R^n,k)\rightarrow B(\R^n,k)$ be the obvious fibration. More generally we can admit $r=\infty$ as well. Then
$$\secat(\pi^n_k)=\cat(B(\R^n,k)\hookrightarrow B(\R^{n+r},k))=\cat_{B(\R^{n+r},k)}B(\R^n,k).$$
\end{thm}

The key to this observation is to consider the ordered Euclidean configuration spaces as well. Another reason for considering the spaces $F(\R^n,k)$ is that $\cat (F(\R^n,k))$ gives a lower bound for $\cat(B(\R^n,k))$ by the usual covering argument. In general this bound is quite bad. For $n=2$ however, it allows to precisely determine $\cat(B(\R^2,k))$ and shows that the subtleties in the calculation for $\secat(\pi^2_k)$ do not arise in the calculation of $\cat(B(\R^2,k))$, which turns out to be $k-1$ for all $k$. The following result might also be of interest for its own sake.

\begin{thm} \label{thmA} 
For all $n\geq 1$, ie, as long as $F(\R ^{n+1},k)$ is connected,
$$\cat(F(\R ^{n+1},k))=k-1.$$ The space $F(\R,k)$ consists of $k!$ contractible
components, hence $$\cat(F(\R,k))=k!-1.$$
\end{thm}

It would be nice to have an analogous statement for the category of the unordered Euclidean configuration spaces, able to compete with the previous theorem in simplicity. We hold the following quite plausible:

\begin{conj} \label{con}
For all $n$ and $k$
$$\cat(B(\R^n,k))=(k-1)\cdot(n-1).$$
\end{conj}
 
Among other indications, our optimism is based on the following calculations. We use the more usual notation $\alpha(k)$ for $D_2(k)$: 

\begin{thm}\label{thmB}
Let $\alpha(k)=D_2(k)$ be the number of {\rm1}'s in the dyadic expansion of $k$.
Then we have
\begin{equation*}
(k-\alpha(k))\cdot(n-1) \leq \secat(\pi^n_k) \leq \cat(B(\R^n,k))
\leq(k-1)\cdot(n-1). \label{general unordered bounds}
\end{equation*}
\end{thm}

\begin{thm}\label{thmC}
In case $k$ is a power of $2$ or $k=3$ or if $n$ is odd and $k=p$ is any prime, we have 
$$\secat(\pi^n_k) = \cat(B(\R^n,k))=(k-1)\cdot(n-1).$$
Moreover, for any $k$ we have $$\cat(B(\R^2,k))=(k-1).$$
\end{thm}

These results are obtained by exploiting work of Vassiliev, cohomology calculations by Fred Cohen and combining them with standard results from LS--theory and the concept of category weight. The statement for $k=3$ follows together with \fullref{sec+sub} and geometric insight. The general upper bound in \fullref{thmB} can be derived from the following lemma, which we could not find in literature.

\begin{lem}\label{myCW}
Let $X$ be an $n$--dimensional CW--complex, $X^{(r)}$ its $r$--skeleton and assume that $X-X^{(k-1)}$ is connected. Then $\cat(X-X^{(k-1)})\leq n-k$.
\end{lem}

We are aware of the incompleteness of \fullref{thmB}. We are also
aware that improvements can be achieved.  For instance, we gained some
generalizations in \cite{Rot05}.  However, we did not obtain a
complete generalization of Vassiliev's results, new improvements of
the upper bound for $\secat(\pi^n_k)$, an unbounded sequence $(n_i)$
such that $\cat(B(\R^{n_i},k))=(k-1)(n_i-1)$ for all $k$ or of course
a proof or disproof of \fullref{con}.


We think it is worth mentioning that the behaviour showing up in Vassiliev's calculations and our \fullref{thmB}, as well as the lack of complete information seem not to be unusual. 
For example, take the immersion problem for real projective spaces $\R P^n$ into $\R^m$ which is still open in the general case. Whitney's embedding theorem says that an immersion exists at least for $m$ greater or equal to $2n-1$. This bound is taken if $n$ is a power of $2$. 
More generally, the number $\alpha(n)$ appears in Ralph Cohen's general immersion theorem which says that every compact, differentiable, $n$--dimensional manifold immerses in Euclidean space of dimension $2n-\alpha(n)$ \cite{Coh85}. 
For the case of complex projective spaces it has been conjectured,
that the immersion dimension is $4n-2\alpha(n)+\epsilon$, where
$\epsilon$ is a non-negative integer bounded roughly by $3$, see
Gonz{\'a}lez \cite{Gon03} and references therein. In fact, these
immersion problems are closely related to invariants of the category
type. For instance, in the real case and $n\not=1,3,7$, Faber,
Tabachnikov and Yuzvinsky showed in \cite{FTY} that the immersion
dimension is the sectional category of $P \R P^n \rightarrow \R P^n
\times \R P^n$, up to shift also known as the topological complexity
of $\R P^n$. Here $P \R P^n$ is the space of all continuous paths
$\gamma \col [0,1] \rightarrow \R P^n$, and the fibration is
evaluation at the end points. This notion was also useful for the
immersion problem for $2^e$--torsion lens-spaces for $e>1$ as an
approach to the immersion problem for complex projective spaces
\cite{Gon03}.

As another example consider the Lusternik--Schirelman category of the real Grassmann manifolds $G_{n,k}$ of $k$ dimensional subspaces in $\R^{n+k}$. By dimensional reasons $\cat(G_{n,k})\leq nk$ and Berstein showed in \cite{Ber76}, that this bound is taken if and only if $n=1$ or $k=1$ or $(n=2$ and $k=2^m-1)$ or $(k=2$ and $n=2^m-1)$. We are not aware of precise determinations in the general case.

This paper developed from the author's diploma thesis \cite{Rot05} which is more detailed in a number of points. However, here we put more emphasis on the sectional category point of view.

\medskip
\begin{fat}
Acknowledgements
\end{fat}\qua
First of all I would like to thank my advisor Carl-Friedrich B\"odigheimer, who suggested to investigate the LS--category of con\-fi\-gu\-ra\-tion spaces and guided my work which lead to \cite{Rot05}.
I am also very grateful to Fred Cohen whom I met first when I started this project and attended his lectures in Louvain-la-Neuve. Besides its mathematical interest, his comments and encouraging remarks were of invaluable mental support. I like to thank the organizers for the wonderful conference and the invitation to Tokyo. I'm also grateful to Daniel Tanr\'e, Sadok Kallel and Yves F\'elix whom I could talk to in Lille about the subject. Jes\'us Gonz\'alez drew my attention to the close relation between sectional category and immersion questions. Special thanks are to Birgit Richter for critical remarks, corrections and discussion of the manuscript. Last but not least I'd like to thank the German National Academic Foundation for their support.

\section{Lusternik--Schnirelman category} \label{LSC}

The investigation of numerical homotopy invariants called \emph{category} began with an article by Lusternik and Schnirelman \cite{LS34}. Their aim was to obtain bounds for the number of critical points of a smooth function on a manifold. Since then various slightly differing definitions showed up in the literature. The definition given in the introduction takes into account that we should have $\cat(*)=0$ for a point or contractible space $*$. We followed \cite{CLOT03} and also recommend this book as a source for the following results. 

We now give a unifying approach to various notions of category including those mentioned above. For that purpose we define the category $\cat(f)$ of a map  $f\col A\rightarrow B$. Let $\cat(f)$ be the least integer $n$, such that $A$ can be covered with $n+1$ open sets and the restriction of $f$ to each of these sets is nullhomotopic. Such a cover of $A$ is called categorical.

We recover the definitions of the introduction via
$\cat(X)=\cat(id_X)$ and also $\secat(p)$ equals $\cat(f_p)$ when
$f_p$ is a classifying map for a principal fibration $p$, see
\cite[9.18, 9.19]{CLOT03} and Hatcher \cite[Exercise 22 page 420]{Hat}. We will
only deal with fibrations of this type and can hence use $\secat(p)$
and $\cat(f_p)$ interchangeably in the sequel. Furthermore, for
$A\subset B$ the subspace category of $\cat_B(A)$ is defined to be the
least $n$ such that there exists a cover of $A$ with $n+1$ subsets of
$B$, each open and contractible in $B$. For an open inclusion $i\col
A\hookrightarrow B$ one obviously has $\cat_B(A)=\cat(i)$.

Alternative definitions for $\cat(X)$ which agree with the standard one under mild hypotheses, including the case where $X$ is a pointed CW--complex, have been given by Whitehead and Ganea. For a space $X$ Ganea constructed a sequence of fibrations $p_n \col G_n(X) \rightarrow X$ which have a section if and only if $\cat(X)\leq n$.

We now recollect some properties of the category of a space $X$:

\begin{prop}\label{LSprop}\ 

\begin{enumerate}
\item
If $X$ dominates $Y$, ie, if there are maps $f\col X\rightarrow Y$ and $g\co Y\rightarrow X$ such that $f\circ g\simeq id_Y$, then $\cat Y\leq \cat X$. In particular, category is a homotopy invariant.

\item \label{cov}
We have $\cat(E)\leq \cat(B)$ for a covering $p\col E\rightarrow B$ with $E$ path-connected.

\item \label{secat}
We have $\secat(p)\leq \cat(B)$ for a fibration $p\col E\rightarrow B$.

\item \label{cupbound}
We have $\cuplength_R(X) \leq \cat(X)$ where $\cuplength_R(X)$ is the $R$--cuplength for any coefficient ring $R$, ie, the least $n$ such that all cup-products of at least n+1 non-trivial factors in $\tilde H^*(X;R)$ vanish.

\item \label{dim}
If $X$ is an $(n-1)$--connected CW--complex, then $\cat(X)\leq \frac{dim(X)}{n}$.
\end{enumerate}
\end{prop}

All proofs are elementary and can be found in \cite[Chapters 1 and
3]{CLOT03}. Since the statement of \fullref{myCW}
fits very well in this collection, we now give its proof:

\begin{proof}
[Proof of \fullref{myCW}]
The idea is the same as one can use to show that every CW--subcomplex is a strong deformation-retract of some open neighborhood. First note that subadditivity (ie, $\cat_X(A\cup B) \leq \cat_X(A)+\cat_X(B) +1$ for $A,B \subset X$) yields 
\begin{eqnarray*}
\cat(X-X^{(k-1)}) &=& \cat_{X-X^{(k-1)}} \left( \coprod_{r=k}^n X^{(r)}-X^{(r-1)} \right) \\
&\leq& n-k + \sum_{r=k}^n \cat_{X-X^{(k-1)}} \left( X^{(r)}-X^{(r-1)} \right).
\end{eqnarray*}
Hence it suffices to show that $\cat_{X-X^{(k-1)}} ( X^{(r)}-X^{(r-1)})=0$ for all $r\in \{k,k+1,\dots,n\}$, ie, that $X^{(r)}-X^{(r-1)}$ is covered by some set which is open and contractible in $X-X^{(k-1)}$. Fix $r\in
\{k,k+1,\dots,n\}$ for the sequel. $X^{(r)}-X^{(r-1)}$ is a disjoint union of $r$--balls, hence contractible in $X-X^{(k-1)}$, since $X-X^{(k-1)}$ is path-connected. We are left to show that $X^{(r)}-X^{(r-1)}$ is a retract of some set $V_n$, open in $X-X^{(k-1)}$. 
The requirement that $V_n$ is contractible in the ambient space $X-X^{(k-1)}$ is then automatically satisfied since $X^{(r)}-X^{(r-1)}$ is a disjoint union of $r$--cells, hence contractible in the ambient space: Each of the cells is contractible to a point and then we use that in our case the ambient space is path-connected.

We will obtain $V_n$ by recursively defining sets $V_{\ell}$ for $r\leq \ell\leq n$, such that $V_{\ell} \supset X^{(r)}-X^{(r-1)}$ is open in $X^{(\ell)}-X^{(k-1)}$ and retracts to $X^{(r)}-X^{(r-1)}$: For $V_r$ we can just take $V_r:=X^{(r)}-X^{(r-1)}$. If $V_{\ell} \supset
X^{(r)}-X^{(r-1)}$ as required is already defined, extend it to obtain $V_{\ell+1}$ as follows: For each $(\ell+1)$--cell $e$ choose a point $x_e$ in its interior, as well as radial homotopies
$h_t^e\co  X^{(\ell)}\cup (e-x_e) \rightarrow X^{(\ell)}\cup (e-x_e)$,
ie, homotopies relative $X^{(\ell)}$ with $h_0^e=id$, $h_1^e\co 
X^{(\ell)}\cup (e-x_e) \rightarrow X^{(\ell)}$ and $h_1^e\circ
h_t^e=h_1^e$. The set $(h_1^e)^{-1}(V_{\ell})$ then is obtained from
$V_{\ell}$ by glueing a truncated cone over $V_{\ell}\cap \partial e$, open in $\overline{e}$. Define $V_{\ell+1}:= \bigcup_e
(h_1^e)^{-1}(V_{\ell})$, taking the union over all $(\ell+1)$--cells $e$.
The set $V_{\ell+1}\supset X^{(r)}-X^{(r-1)}$ then is open in $X^{(\ell+1)}-X^{(k-1)}$ (weak topology) and retractible to $X^{(r)}-X^{(r-1)}$. This shows $\cat_{X-X^{(k-1)}} \left(
X^{(r)}-X^{(r-1)} \right)=0$. 
\end{proof}

Most of our lower bounds for the sectional category are obtained
through the concept of category weight $\wgt_R$. For a non-zero class
$u\in H^*(X;R)$ define $\wgt_R(u)$ to be the greatest $k$ (or
$\infty$), such that $p_{k-1}^*(u)=0 \in H^*(G_{k-1}(X);R)$ for the
$(k-1)$st Ganea fibration $p_{k-1}\co G_{k-1}(X)\rightarrow X$. We
recollect important properties and consequences from \cite[pages 63f,
242ff and 261f]{CLOT03}. The last point is a consequence of
\cite[Proposition  9.18 and 8.22(2)]{CLOT03}:

\begin{prop}\label{wgt}
Let $u\in H^k(X;R)$ be non-zero. Then:
\begin{enumerate}

\item If $u\in H^k(K(\pi,1);R)$ 
is a class in the cohomology of an Eilenberg--MacLane space of type $(\pi,1)$, then $\wgt_R(u)= k$. \label {KP1}
\item If $f\co Y\rightarrow X$ is such that $f^*(u)\not =0$, then $\wgt_R (f^*(u))\geq \wgt_R(u)$.
In other words: If a cohomology class does not vanish under pullback, then its category weight cannot decrease. \label{rueck}
\item If $p\col E\rightarrow Y$ is a fibration arising as a pullback over $f\col Y \rightarrow X$ of a fibration $\hat E \rightarrow X$ with contractible total space $\hat E$ and $f^*(u)\not =0$, then $\wgt_R(u) \leq \secat (p)$. \label{absch}
\end{enumerate}
\end{prop}

\section{Cellular models, geometry and cohomology of Euclidean configuration spaces}

In this section we will collect the necessary algebraic and geometric data in order to derive bounds for the category of Euclidean configuration spaces in combination with the results of \fullref{LSC}. The maps and spaces under consideration fit into the following fundamental diagram

\begin{equation}\label{funddiag}
\xymatrix{\sym_k \ar[d] & \sym_k \ar[d] && \sym_k \ar[d] \\
                F(\R^n,k) \ar[r]^(0.47){\tilde f} \ar[d]_{\pi^n_k} & F(\R^{n+1},k)\ar[r] \ar[d]^{\pi^{n+1}_k} &\cdots \ar[r]  &F(\R^{\infty},k)\simeq * \ar[d] ^{\pi^{\infty}_k}\\
                B(\R^n,k) \ar[r]^(0.45)f & B(\R^{n+1},k)\ar[r] &\cdots \ar[r]& B(\R^{\infty},k)=K(\sym_k,1) }
\end{equation}
 
where the vertical maps are the coverings given by the free action of the symmetric group $\sym_k$ and horizontal maps are induced by the inclusion of $\R^n$ into $\R^{n+1}$. It follows from Fadell's and Neuwirth's fundamental sequence of fibrations \cite{FN62}, that $F(\R^n,k)$ is $n-2$ connected. As a consequence the limit spaces on the right give the universal covering of an Eilenberg--MacLane space $K(\sym_k,1)$. Note for later use that all the rectangles in diagram \eqref{funddiag} are homotopy pullbacks.

The integer cohomology of Euclidean configuration spaces was calculated by Fred Cohen \cite{CLM76,CT78}. In the formulation of \cite{Coh95} and for $n\geq2$, $H^*(F(\R^n,k))$ is given by generators $A_{i,j}$, $(1\leq j<i \leq k)$ all in degree $n-1$, subject to the relations 
\begin{enumerate}
\item $A_{i,j}^2=0$
\item $A_{i,j}A_{i,\ell}=A_{\ell,j}(A_{i,\ell}-A_{i,j})$  for $j<\ell<i$
\item associativity and graded commutativity.
\end{enumerate}

We can draw the following conclusions:

\begin{cor} \label{cup}
For all $n\geq2$, we have $\cuplength_{\Z}(F(\R^n,k))=k-1$.
\end{cor}

\begin{cor}
Let $n\geq 3$. Then $F(\R^n,k)$ is homotopy equivalent to a CW--complex $Y$ with cells just in dimension $q \cdot (n-1)$ for $q\in \{0,1,\dots,k-1\}$.
\end{cor}

\begin{proof}
Follow the construction in \cite[4.C page 429]{Hat} and note that in our case the homology is free. A geometric construction (also for $n=2$) can also be found in \cite[Sections VI.8,VI.10]{FH01}. 
\end{proof}

These corollaries will allow to calculate the category of ordered configuration spaces completely and we now turn to the unordered case and sectional category. 

We are going to apply the concept of category weight using the fundamental diagram \eqref{funddiag}. 
It turns out that we can draw a lot of information from a CW--decomposition of the one-point compactification $B(\R^n,k)_{\infty}$ of $B(\R^n,k)$ introduced by Vassiliev \cite{Vas88}\cite[pages 28ff]{Vas92}. This decomposition is a generalization of the one introduced by Fuks for $n=2$ in \cite{Fuk70}. Vassiliev describes the various cells of his model as well as their boundaries mod $2$ in terms of certain so-called $(n,k)$--trees. A precise description of what an $(n,k)$--tree looks like is given in the construction \cite[page 28]{Vas92}. We only note that an $(n,k)$--tree has at least $k+n-1$ and at most $k\cdot n$ edges. 
Vassiliev has proven that for any $n,k$ there exists the structure of a $CW$--complex of the space $B(\R^n,k)_{\infty}$ with cells being sets of points corresponding to various $(n,k)$--trees and the added point \cite[lemma 3.3.1, page 29]{Vas92}. Furthermore, the dimension of such a cell is equal to the number of edges in the corresponding tree \cite[lemma 3.3.2, page 29]{Vas92}. Altogether, this leads to the following observation:

\begin{prop}[Vassiliev] \label{VasCW}
There is a CW--decomposition of the one-point compacti\-fi\-ca\-tion $B(\R^n,k)_{\infty}$ having the point $\infty$ as the only cell of dimension $0$. All other cells have dimension $r$ with $k+n-1\leq r \leq k\cdot n$.
\end{prop}

There is a stabilization of these models as $n$ turns to $\infty$, and using Poincar\'e--Lefschetz duality, Vassiliev shows:

\begin{prop}{\rm \cite[page 27]{Vas92}}\label{sur}\qua
The homomorphism $$H^*(B(\R^{\infty},k);\zz)\rightarrow H^*(B(\R^n,k);\zz)$$ induced by the map from the fundamental diagram \eqref{funddiag} is surjective.
\end{prop}

The reader who is familiar with Vassiliev's cell decomposition and its description in terms of trees may also derive the following:

\begin{cor} \label {cohdim}
Let $\alpha(k)$ be the number of $1$'s in the dyadic decomposition of $k$. Then
$$H^q(B(\R^n,k);\zz) \left\{ \begin{array}{ll} =0 & \quad \mbox{if $q>(k-\alpha(k))\cdot(n-1)$}\\
                                                                  \not =0 &\quad \mbox{if $q=(k-\alpha(k))\cdot(n-1)$.}
                                                                            \end{array} \right.$$
In other words
$$ \cohdim_{\zz}B(\R^n,k)=(k-\alpha(k))\cdot(n-1).$$
\end{cor}

\begin{proof}
For $k=2$ this follows from \cite[Section 4,4.2]{Fuk70} with an
elementary proof on page 144f. For arbitrary $k\geq 2$ this follows
from the theorem in \cite[page 31]{Vas92} once one is familiar with
Vassiliev's cell decomposition. We do not want to repeat this
construction but give some hints for the reader who wants to get
acquainted with Vassiliev's notation and the labeling trees
$\Gamma_{K_i}$ that occur in his theorem: For the case of
$H^*(B(\R^n,k))$ (which is $H^*(\R^n(k))=H^*(\R^n(m))$ in Vassiliev's
notation, hence $k=m$) such a tree has vertices concentrated on $n+1$
horizontal lines and branches from top to bottom. A typical example is
given in \cite[Figure 12, page 31]{Vas92}.
The depth of such a tree is the number of the highest horizontal line beneath which no more branchings exist, see \cite[Figure 11, page 30]{Vas92}. The branching condition (page 30, bottom) implies that the number of vertices on each horizontal line is a power of two and the tree $\Gamma_{K_i}$ has $2^{|K_i|}$ vertices on the bottom line. The $m$ in $[K_1,\dots, K_l;m]$ stands for $m-2^{|K_1|}-\cdots -2^{|K_l|}$ copies of the unique tree $\Gamma_0$ that is just a vertical chain of edges without any branching. This means that if we add up the number of vertices on the bottom horizontal lines over all the trees of the collection $[K_1,\dots, K_l;m]$, the sum is $m$. We are looking for an additive generator of maximal degree, which means that we are looking for a collection $[K_1,\dots, K_l;m]$ where the sum of the edges of all the trees is minimal. This is because $H^*(B(\R^n,k))$ is obtained via Poincar\'e duality from the space $B(\R^n,k)_{\infty}$ whose cells correspond to trees, its dimensions correspond to the number of edges. 
From the fact that on each bottom line the number of vertices is a power of two and their sum is $k=m$, it follows that a collection $[K_1,\dots, K_l;m]$ consists of at least $\alpha(m)=\alpha(k)$ trees. For $k=2^{l_1}+\cdots +2^{l_{\alpha(k)}}$ such a collection of trees with the minimal number of edges is given in \fullref{tree}.

\begin{figure}[ht!]
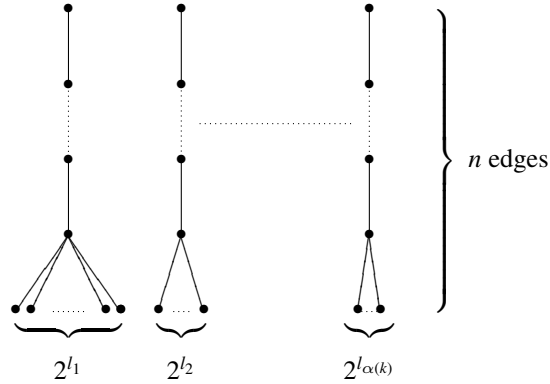
\small
$$\xy 0;<.5cm,0cm>:
@={(3,8),(3,6),(3,4),(3,2),(0,8),(0,6),(0,4),(0,2),(8,8),(8,6),(8,4),(8,2),(7.7,0),(8.3,0),
(2.4,0),(3.6,0),
(-1,0),(-1.4,0),(1,0),(1.4,0),},
@@{*{\bullet}};
(0,2); (0,4)**@{-};
(0,6)**@{.};
(0,8)**@{-};
(3,8)
;
(3,6)**@{-};
(3,4)**@{.};
(3,2)**@{-};
(7.7,0);(8,2)**@{-};
(8,4)**@{-};
(8,6)**@{.};
(8,8)**@{-};
(8.3,0);(8,2)**@{-};
(-1.0,0);(0,2)**@{-};
(-1.4,0)**@{-};
(1,0);(0,2)**@{-};
(1.4,0)**@{-};
(2.4,0);(3,2)**@{-};
(3.6,0);(3,2)**@{-};
%
%
(-.4,0);(.4,0)**@{.};
(2.8,0);(3.2,0)**@{.};
(7.8,0);(8.2,0)**@{.};
(3.5,5);(7.5,5)**@{.};
(0,0);(10,8)
!<0cm,-2.0cm>;
**\frm{\}};
(11,4) *\txt{\qquad \mbox{$n$ edges}};
(-1.4,-.5);(1.4,-.5);
!<.7cm,0cm>;
**\frm{_\}};
(0,-1.5) *\txt{$2^{l_1}$};
(2.4,-.5);(3.6,-.5);
!<.3cm,0cm>;
**\frm{_\}};
(3,-1.5) *\txt{$2^{l_2}$};
(7.8,-.5);(8.2,-.5);
!<.1cm,0cm>;
**\frm{_\}};
(8,-1.5) *\txt{$2^{l_{\alpha(k)}}$};

\endxy$$
\caption{\label{tree}A collection of trees representing an additive generator in\newline $H^{(n-1)(k-\alpha(k))}(B(\R^n,k),\zz)$.}
\end{figure}

This collection has $k+(n-1)\cdot \alpha(k)$ edges and hence represents a generator of $H^*(B(\R^n,k))$ in degree $nk-(k+(n-1)\cdot \alpha(k))=(n-1)\cdot (k-\alpha(k))$. 
More details, pictures and examples are also given in \cite[pages 25--36]{Rot05}. In particular, see Satz 2.9 on page 34.
\end{proof}

There is also a statement for $\zp$--coefficients available in
literature, which we shall exploit later. It goes back to Fred Cohen
\cite{Cohen1,CLM76} and was proved anew and stated explicitly by Ossa
(see \cite[Proposition 3.4]{Oss96} and the following remark for
$p=3$):

\begin{prop} \label{ossassatz}
If $p$ is an odd prime, then
$$H^{(p-1)(n-1)}(\sym_p;\zp)\longrightarrow
H^{(p-1)(n-1)}(B(\R^n,p);\zp)$$ is an isomorphism.
\end{prop}

The reader interested in cohomological dimensions of configuration
spaces should have a look at Kallel \cite{Kal}.

\section{Calculations and proofs}

We begin with the partly special arguments for the computation of $\cat(F(\R^n,k))$ and $\cat(B(\R^n,k))$ in case $n$ or $k$ is less or equal to $2$. 

\begin{fat}$k=1$\end{fat}\qua
Clearly, $\cat(F(\R^n,1))=\cat (B(\R^n,1))=\cat(\R^n)=0$. 

\begin{fat}$n=1$\end{fat}\qua
If $n=1$, then $F(\R,k)$ has $k!$ contractible components and $B(\R,k)$ is contractible, hence $\cat(F(\R,k))=k!-1$ and $\cat (B(\R,k))=0$. 

\begin{fat}$k=2$\end{fat}\qua
Next, since $\R^n$ is a topological group we have $F(\R^n,2)\cong \mathbb R^n \times F(\mathbb R^n-\{0\},1)\simeq \mathbb S^{n-1}$ \cite{CLM76} and hence $\cat (F(\mathbb R^n,2)) = 1$ for all $n\geq 1$.
In the unordered case we use the homotopy equivalence $B(\R^n,2) \simeq \R P^{n-1}$ \cite{CLM76} which implies $\cat (B(\R^n,2))=\cat(\R P^{n-1})=n-1$, bounded above by the dimension and below by the mod $2$ cuplength.

\begin{fat}$n=2$\end{fat}\qua
We have $k-1=\cuplength_{\mathbb Z}(F(\R^2,k))\leq
\cat(F(\R^2,k))\leq \cat (B(\R^2,k))\leq k-1$ resulting from \fullref{cup}, \fullref{LSprop} \eqref{cupbound} and \eqref{cov} and the more general \fullref{upperbound}. 

\begin{lem}\label{upperbound}
For all $n$ and $k$ we have $\cat (B(\R^n,k))\leq (k-1)\cdot (n-1)$.
\end{lem}

\begin{proof}
We can assume $n \geq 2$. Then the lemma follows from \fullref{VasCW} and \fullref{myCW}.
\end{proof}

\begin{proof}[Proof of \fullref{thmA}]
We now can assume $n,k\geq 3$. Then $k-1=\cuplength_{\mathbb Z}(F(\R^n,k))\leq \cat(F(\R^n,k))\leq k-1$ follows from \ref{cup}, \ref{LSprop}\eqref{cupbound} and \eqref{dim}. Remember that $F(\R^n,k)$ is $(n-2)$--connected. 
\end{proof}

\begin{proof}[Proof of \fullref{thmB}]
It follows from \fullref{sur} together with \fullref{wgt}, that the degree of each non-zero cohomology class in $H^*(B(\R^n,k);\zz)$ gives a lower bound for $\secat(\pi^n_k)$. Hence by \fullref{cohdim} and the previous \fullref{upperbound} we have $(k-\alpha(k))\cdot(n-1)=\cohdim_{\zz}(B(\R^n,k)) \leq \secat(\pi^n_k) \leq \cat(B(\R^n,k)) \leq(k-1)\cdot(n-1)$. Here we use \ref{LSprop}\eqref{secat} and \fullref{upperbound} again. 
\end{proof}

\begin{proof}[Proof of \fullref{sec+sub}]
For notational convenience we just formulate the proof for $r=1$. The other cases are similar. We first show $\secat(\pi^n_k)=\cat(B(\R^n,k)\hookrightarrow B(\R^{n+1},k))$ and denote the inclusion by $f$ as in the fundamental diagram \eqref{funddiag}. 
Consider $\secat(\pi^n_k)$ as the category of a classifying map $B(\R^n,k)\rightarrow B(\R^{\infty},k)$. 
The classifying map of $\pi^n_k$ factors through $f$, hence $\secat(\pi^n_k)\leq \cat(f)$. 

Now consider the left square in diagram \eqref{funddiag}. Given a subset $A\subset B(\R^n,k)$ over which $\pi^n_k$ is trivial, we can factor $f|_A$ as $\pi^{n+1}_{k} \circ \tilde f \circ s$ where $s$ is a local section for $\pi^n_k$ over $A$. Now we observe that $F(\R^n,k)$ is contractible in $F(\R^{n+1},k)$: 
Remember that a point in $F(\R^n,k)$ is a $k$--tuple in $\R^n\cong \R^n\times\{0\} \subset \R^{n+1}$. First move the $k$ points of such a tuple linearly by varying only their last coordinates such that in the end the $i^{\tth}$ point lies in $\R^n \times\{i\}$. Then continue moving the $i^{\tth}$ point linearly to $(0,\dots,0,i)$. 
The fact that $F(\R^n,k)$ is contractible in $F(\R^{n+1},k)$ implies that the restriction of $f$ to $A$ is nullhomotopic. This shows $\cat(f)\leq\secat(\pi^n_k)$.

The second equality in \fullref{sec+sub} can be proved quite elementary by pulling back and extending categorical covers. Remember that a point in $U_i\subset B(\R^n,k)$ is a $k$--element subset of $\R^n\cong \R^n\times \{0\}$. If $U_i$ is open and contractible, extend it to an open contractible subset of $B(\R^{n+1},k)$ by simply letting vary the elements of its points in their $(n+1)^{\tst}$ coordinates within an open interval, say $(-1,1)$.
Alternatively one can apply the next lemma and the fact that the spaces $B(\R^n,k)$ are normal (as CW--complexes) and absolute neighborhood retracts (ANR's) as retracts of some open subset of some $\R^N$. For a definition of ANR see the appendix of \cite{CLOT03} and Warner \cite{War} for a more detailed introduction. 
\end{proof}


%

\begin{lem}\label{catcl}
If $i\col A\hookrightarrow B$ is a closed inclusion between normal ANR's, then 
$\cat_B(A)=\cat(i).$
\end{lem}

\begin{proof}
Given a categorical cover for $\cat_B(A)$, inverse images under $i$ give a categorical cover for $\cat(i)$, hence $\cat(i)\leq \cat_B(A)$. Vice versa, given a categorical cover $U_0,\cdots,U_k$ for $\cat(i)$, we can pass to an open refinement $V_0,\cdots,V_k$ with $V_i\subset \overline V_i \subset U_i$ since $A$ is normal \cite[Theorem A.1]{CLOT03}. Hence $A$ can be covered with $k+1$ sets, each closed and contractible in $B$. In \cite{CLOT03} this fact is denoted by $\catcl_B(A)\leq k$ and under the hypothesis that $B$ is a normal ANR and $A\subset B$ is closed,  \cite[Theorem 1.10 ]{CLOT03} says that $\catcl_B(A)=\cat_B(A)$. Hence we have $\cat_B(A)\leq \cat(i)$. 
\end{proof}

\begin{proof}[Proof of \fullref{thmC}]
The statement for $k$ a power of $2$ is a corollary to \fullref{thmB}. If $k=p$ is an odd prime, we combine \fullref{ossassatz} with the group cohomology $H^*(\sym_p;\zp)\cong
\zp[\alpha] \otimes \Lambda_p(\beta)$ where $\alpha$ is a polynomial
generator in degree $2(p-1)$ and $\beta$ is an exterior generator in
degree $2(p-1)-1$. This can be derived from Adem and Milgram
\cite[VI.1.4,1.6, III.2.9]{AM94} or see Ossa \cite{Oss96} for the
statement. Hence $H^{(p-1)(n-1)}(\sym_p;\zp)\neq 0$ if $n$ is odd. Now
we can argue as in the proof of \fullref{thmB} and obtain
$$\secat(\pi^n_p)\geq \left\{ \begin{array}{ll} (p-1)\cdot(n-1) &
\quad \mbox{if $n$ is odd}\\ (p-1)\cdot(n-2) &\quad \mbox{if $n$ is
even.}  \end{array} \right.$$ The inequality for even $n$ is a
consequence of the statement for odd $n$, since
$$H^*(\sym_p;\zp)\rightarrow H^*(B(\R^{n-1},p);\zp)$$ factors via
$H^*(B(\mathbb R^n,p);\zp)$.  The improvement for $k=3$ then follows
from the next lemma and the statement $\cat(B(\R^2,k))=(k-1)$ was
already shown at the beginning of this section.
\end{proof}

\begin{lem}\quad \label{m=3}
For all $n$ we have
\begin{equation}
\secat(\pi^{n+1}_3)\leq \secat(\pi^n_3) +2.
\end{equation}
\end{lem}

\begin{proof}
We partition $B(\R^{n+1},3)=\coprod_{k=1}^3 V_k(\R^{n+1})$, where
$V_k(\R^{n+1})$ is the sub\-manifold of all $3$--configurations in $\R^{n+1}$ whose image under the perpendicular projection onto $\R^n\cong \R^n \times \{0\} \subset\R^{n+1}$ consists of exactly $k$ points. 
$V_1(\R^{n+1})$ is obiously contractible. 
The space $V_2(\R^{n+1})$ is not necessarily contractible but it is contractible whithin $B(\R^{n+1},3)$. This can be seen by an argument similar to the one that we used in order to show that $F(\R^n,k)$ is contractible within $F(\R^{n+1},k)$. The space $V_2(\R^n)$ can be contracted within $B(\R^{n+1},3)$ in the following way: Move the three points (making up a point in $V_2(\R^n)$) linearly by varying only their last coordinates to obtain a three-element subset of $\R^{n+1}$ of the form $\{(x,-1),(\tilde x,0),(\tilde x,1)\}$ (here $x,\tilde x \in \R^n$), then move this linearly to $\{(0,0,-1),(0,0,0),(0,0,1)\}$. Furthermore $\cat_{B(\R^{n+1},3)}(V_3(\R^{n+1}))\leq \cat (B(\R^n,3))$, since $V_3(\R^{n+1})\subset B(\R^{n+1},3)$ is open and retractible to $B(\R^n,3)$. Now we should pass to tubular neighborhoods $U_1, U_2$ of $V_1,V_2$ in order to have \emph{open} contractible sets available. We obtain 
\begin{eqnarray*}
\secat(\pi^{n+1}_3) &\leq& \cat(B(\R^{n+1},3)) \\
&=& \cat(\coprod_{k=1}^3 V_k(\R^{n+1})) \\
&=&\cat_{B(\R^{n+1},3)}(U_1\cup U_2 \cup V_3) \\
&\leq& \cat_{B(\R^{n+1},3)}(V_3)+\cat_{B(\R^{n+1},3)}(U_2)+\cat_{B(\R^{n+1},3)}(U_1)+2 \\
&\leq& \cat_{B(\R^{n+1},3)}(B(\R^{n},3)) +2 \\
&=& \secat(\pi^n_3) +2.
\end{eqnarray*}
The last inequality follows from \fullref{sec+sub} with $r=1$. 
\end{proof}
\bibliographystyle{gtart}
\bibliography{link}

\end{document}